\documentclass[preprint]{sig-alternate}

\usepackage{dsfont}

\newtheorem{Cond}{Condition}
\newtheorem{far}[Cond]{Theorem}
\newtheorem{mn}[Cond]{Corollary}
\newtheorem{clus}[Cond]{Theorem}
\newtheorem{rm1}[Cond]{Remark}

\newtheorem{lem1}[Cond]{Lemma}
\newtheorem{lem2}[Cond]{Lemma}
\newtheorem{chen}[Cond]{Theorem}

\begin{document}
%
\CopyrightYear{} 
\crdata{}  

\title{Far-out Vertices In Weighted Repeated Configuration Model}
%
%
%
%
%

\numberofauthors{2} 
%
\author{
%
%
\alignauthor
B.~B{\l}aszczyszyn
\\
       \affaddr{INRIA/ENS}\\
       \affaddr{23 avenue d'Italie}\\
       \affaddr{Paris, France}\\
       \email{Bartek.Blaszczyszyn@ens.fr}
\alignauthor
K. Gaurav
\\
       \affaddr{UPMC/INRIA/ENS}\\
       \affaddr{23 avenue d'Italie}\\
       \affaddr{Paris, France}\\
       \email{kumar.gaurav@inria.fr}
}
\date{30 July 1999}

\toappear{Appeared in the proceedings of the MAMA workshop, held in conjunction with ACM Sigmetrics/Performance 2012.}

\maketitle
\begin{abstract}
We consider an edge-weighted uniform random graph with a given degree
sequence (Repeated Configuration Model) which is a useful approximation for many real-world networks.
It has been observed that the vertices which are separated from the rest of the graph by a distance exceeding certain threshold play an important role in determining some global properties of the graph like diameter,
flooding time etc., in spite of being statistically rare. We give a convergence result for the distribution of the number of such far-out vertices. We also make a conjecture about how this relates to the
longest edge of the minimal spanning tree on the graph under consideration.
\end{abstract}




\section{Introduction}
Given a degree sequence $\text{\bf d}^{(n)}=(d_i^{(n)})_1^n$ for $n$ vertices labelled $1$ to $n$, Repeated Configuration Model, denoted $RCM[n]$ (or $RCM[(\text{\bf d}^{(n)})]$ when talking about multiple models each with different degree sequence), gives a uniform random graph with vertex $i$ having degree $d_i$. 
Every edge $e$ of $RCM[n]$ is given an i.i.d. edge-length $Y(e)$, with cdf $F$, to get a weighted random graph which we denote by $\widetilde{RCM}[n]$.
Given $x_n>0$, a vertex $i$ in $\widetilde{RCM}[n]$ is said to be far-out if the distance of $i$ from its nearest neighbour exceeds $x_n$. That is, if we let
$ M_n(i) := \min_{j:j\sim i}Y(e_{ij})$, $i$ is far-out if $M_n(i)>x_n$. In this work, we study the asymptotic distribution of the number of far-out vertices in $\widetilde{RCM}[n]$ as $x_n$ is appropiately scaled with $n$.
This study is motivated by the following observations:
\begin{itemize}
\renewcommand{\labelitemi}{$\bullet$}
\item In the study of the diameter of $\widetilde{RCM}[n]$ in~\cite{Marc}, the existence of far-out vertices plays a critical role in determining the difference between the expected distance between two uniformly chosen vertices and the diameter of the graph.
\item Penrose shows in~\cite{Penrose} that for $n$ points placed uniformly at random in the unit square, the longest edge of \textit{Minimal Spanning Tree(MST)} and the longest edge of \textit{Nearest Neighbour Graph(NNG)} on these points are asymptotically the same.
In other words, the most far-out vertex in the complete graph on these $n$ points determines the longest edge of the MST of this complete graph.
\end{itemize}
 
The second observation suggests a conjecture  that the longest edge of the MST on $\widetilde{RCM}[n]$ would be determined by the most far-out vertex in $\widetilde{RCM}[n]$.
However in the present abstract, we limit ourselves to presenting the results on the asymptotoic distribution of the number of far-out vertices and the asymptotic distribution of the longest
edge of the NNG of $\widetilde{RCM}[n]$.

\section{Results}
For a precise statement of the results, we start by assuming a set of conditions on $\text{\bf d}^{(n)}=(d_i^{(n)})_1^n$ given in chapter 7 of~\cite{Remco}, which apart from being indispensable for convergence results, enables the construction of
$\widetilde{RCM}[n]$ for arbitrary large $n$. 
\begin{Cond}
\label{degden}
For each $n$, $\text{\bf d}^{(n)}=(d_i^{(n)})_1^n$ is a sequence of non-negative integers such that $\sum_{i=1}^nd_i^{(n)}$ is even. For $k \in \mathbb{N}$,let $u_k^{(n)}=\vert \{i:d_i^{(n)}=k\} \vert$ and $D_n$ be the
degree of a uniformly chosen vertex in $\widetilde{RCM}[n]$, i.e., $\mathbb{P}(D_n=k)=u_k^{(n)}/n$. Let $D$
be some random variable taking value over strictly positive integers with probability distribution $(p_r)_1^{\infty}$.
Then the following hold.
\begin{enumerate}
 \item The degree density condition: $u_k^{(n)}/n \to p_k$ for every $k$ or equivalently, $D_n \to D$ in distribution.
 \item Convergence of emperical expectations:  $\sum_{i=1}^n \frac{d_i^{(n)}}{n} \to  \mathbb{E}[D] \in (0,\infty)$.
\end{enumerate}
\end{Cond}


Let $d_{min}:= \min \{d:\mathbb{P}(D=d)>0\}$, $\bar{F}=1-F$ and 
 $N_n$ be the number of far-out vertices in $\widetilde{RCM}[n]$. In what follows, we assume that $d_{min} \geq 3$, $d_i^{(n)} \geq d_{min}$ for all $i$ and $n$ and that $\bar{F}$ is continuous. We have then,

\begin{far}
\label{farout}
 Suppose the sequence $(x_n)_{n\in \mathbb{N}}$ is such that when $n\to \infty$,
 \begin{equation}
 \label{xncon}
  n\bar{F}(x_n)^{d_{min}} \to x
 \end{equation}
for some fixed
$x > 0$. Then $N_n$ converges in the variation distance topology to the Poisson random variable
with parameter $xp_{d_{min}}$ for almost all sequences $\{RCM[n]\}_{n \geq 1} $.
\end{far}

\begin{mn}
\label{mncor}
Let $\text{\bf M}_n:=\max_{1\leq i \leq n}M_n(i)$. Then, for $x>0$
\begin{equation}
\label{mnex}
\mathbb{P}(n\bar{F}(\text{\bf M}_n)^{d_{min}} > x) \to e^{-xp_{d_{min}}}
\end{equation}

\end{mn}

Heuristically, given a particular realization of $RCM[n]$, $I_i^{(n)}:=\mathds{1}(\text{vertex } i \text{ is far-out})$ is a Bernoulli random variable with parameter $\bar{F}(x_n)^{d_i^{(n)}}$ and
$N_n=\sum_{i=1}^n I_i^{(n)}$. Ignoring the weak dependence among $I_i^{(n)}$'s,
\begin{equation*}
N_n \stackrel{d}{\approx} Binomial(n,\sum_{i=1}^n \bar{F}(x_n)^{d_i^n})
\end{equation*}

Now,
\begin{equation}
\label{dec}
 \sum_{i=1}^n \bar{F}(x_n)^{d_i^n} \approx n\mathbb{E}[\bar{F}(x_n)^D] =  n\bar{F}(x_n)^{d_{min}}\mathbb{E}[\bar{F}(x_n)^{\tilde{D}}]
\end{equation}
where $\tilde{D}=D-d_{min}$.
From (\ref{xncon}), when $n \to \infty$, $\bar{F}(x_n) \to 0$ and therefore
\begin{equation*}
  \mathbb{E}[\bar{F}(x_n)^{\tilde{D}}] \to \mathbb{P}(\tilde{D}=0) = \mathbb{P}(D=d_{min}) = p_{d_{min}}
\end{equation*}
This, along with (\ref{xncon}) and (\ref{dec}), implies that for large $n$,\\ 
$\sum_{i=1}^n \bar{F}(x_n)^{d_i^n} \approx xp_{d_{min}}$  and hence,
\begin{equation*}
 N_n \stackrel{d}{\approx} Poisson(xp_{d_{min}})
\end{equation*}

We make the above steps rigorous using the Stein-Chen method for which we need the following theorem and lemmae. Here, $P_a$ denotes Poission r.v. with parameter $a$, $\mathcal{L}(\boldsymbol{V})$ denotes the law of 
random vector $\boldsymbol{V}$, and $d_{var}(.,.)$ denotes the \textit{distance in total variation}. 

\begin{chen}[Stein-Chen]
 For a finite or countable set $V$, let $(I_j)_{j\in V}$ be a family of not necessarily
independent Bernoulli variables with respective parameters $({\pi}_j)_{j\in V}$.
Let $X=\sum_{j\in V} I_j$ and $\lambda=\sum_{j\in V}\pi_j$.
Assume that there exist random variables $(J_{ij})_{i,j \in V, j\neq i}$ defined
on the same probability space as $(I_j)_{j\in V}$ and such that for all $i \in V$,
the following equality of distributions holds:
\begin{equation}
\label{chenloi}
 \mathcal{L}((J_{ij})_{i,j \in V, j\neq i})=\mathcal{L}((I_j)_{j\in V, j\neq i}|I_i=1)
\end{equation}
 Then
\begin{equation}
\label{sten}
 d_{var}(\mathcal{L}(X),P_{\lambda})\leq 2 \frac{1-e^{-\lambda}}{\lambda} \sum_{i \in V}\pi_i(\pi_i + \sum_{j\in V, j\neq i}\mathbb{E}(|I_j-J_{ij}|)
\end{equation}

\end{chen}
\proof
see chapter 2,~\cite{mass}
\begin{lem1}
\label{lemma1}
 For $\lambda,\mu \geq 0$,
\begin{equation*}
 d_{var}(P_{\lambda},P_{\mu})\leq 2|\lambda-\mu|
\end{equation*}

\end{lem1}
\proof
see chapter 2,~\cite{mass}

\begin{lem2}
\label{lemma2} 
Given the Condition~\ref{degden} and the condition on $(x_n)_{n\in \mathbb{N}}$ in Theorem~\ref{farout}, we have that when $n\to \infty$
 \begin{equation}
   \sum_i \bar{F}(x_n)^{2d_i^{(n)}} \to 0
 \end{equation}

\begin{equation}
 \sum_i d_i^{(n)} \bar{F}(x_n)^{d_i^{(n)}+1} \to 0
\end{equation}

\end{lem2}

\proof
We have from (\ref{xncon}), when $n \to \infty$, 
\begin{equation}
\label{dicon}
 \bar{F}(x_n) \to 0 
\end{equation}
Then,
\begin{align*}
\sum_i \bar{F}(x_n)^{2d_i^{(n)}} &\leq n \bar{F}(x_n)^{2d_{min}} \\
&\leq \bar{F}(x_n)^2 n \bar{F}(x_n)^{d_{min}} \\
\end{align*}
which converges to $0$ by (\ref{xncon}) and (\ref{dicon}).

Similarly,
\begin{equation*}
\sum_i d_i^{(n)} \bar{F}(x_n)^{d_i^{(n)}+1} \leq \bar{F}(x_n) n \bar{F}(x_n)^{d_{min}} \sum_i \frac{d_i^{(n)}}{n} 
\end{equation*}
which converges to $0$ by (\ref{xncon}), (\ref{dicon}) and Condition~\ref{degden}.

\vspace{5mm} 

\proof[of theorem~\ref{farout}]
 
 \quad \\
 \quad Let $I_i^{(n)}:=\mathds{1}(\text{vertex } i \text{ is far-out in } \widetilde{RCM}[n])$ be a Bernoulli random variable with parameter $\pi_i:=\bar{F}(x_n)^{d_i^{(n)}}$. Therefore,
$N_n=\sum_{i=1}^n I_i^{(n)}$ and let $S = \sum_i \bar{F}(x_n)^{d_i^{(n)}}$. We will drop the supercript $(n)$, where it is clear from context. Further, let
\begin{equation*}
 J_{ij}:=\mathds{1}(\min_{k:k\sim j,k\neq i}Y(e_{kj}) > x_n)
\end{equation*}
From the independence of edge weights, it is clear that (\ref{chenloi}) is satisfied, given $RCM[n]$.
Now,
\begin{multline*}
 |I_j-J_{ij}|=J_{ij}-I_j=\mathds{1}(\exists \text{ edge } e_{ij} \text{ in } RCM[n], Y(e_{ij}) \leq x_n \\
\text{ and } \min_{k:k\sim j,k\neq i}Y(e_{kj}) > x_n)
\end{multline*}
Therefore,
\begin{align*}
 \mathbb{E}[|I_j-J_{ij}|\big|RCM[n]] = \mathds{1}(\exists \text{ } e_{ij}\text{ in } RCM[n])F(x_n)\bar{F}(x_n)^{d_j-1}
\end{align*}
Plugging this into (\ref{sten}), we get
\begin{align*}
  d_{var}(\mathcal{L}(N_n\big|RCM[n]),P_{S}) &\leq 2 \frac{1-e^{-S}}{S}[\sum_i \pi_i^2 \\
&\quad+ \sum_i \pi_i\sum_{j\neq i}\mathds{1}(\exists \text{ } e_{ij})F(x_n)\bar{F}(x_n)^{d_j-1}] \\ 
&\leq 2 \min(1,\frac{1}{S})[\sum_i\bar{F}(x_n)^{2d_i}\\
&\quad+\sum_i \bar{F}(x_n)^{d_i} \sum_{j\sim i}F(x_n)\bar{F}(x_n)^{d_j-1}] \\
&\leq 2 \min(1,\frac{1}{S})[\sum_i\bar{F}(x_n)^{2d_i}\\
&\quad+\sum_i \bar{F}(x_n)^{d_i+1} d_i]
\end{align*}
Thus, by Lemma~\ref{lemma2}, we have when $n\to \infty$ 
\begin{equation}
\label{dv1}
  d_{var}(\mathcal{L}(N_n\big|RCM[n]),P_{S}) \to 0
\end{equation}

Now, 
\begin{align*}
\label{dec}
 S &=\sum_{i=1}^n \bar{F}(x_n)^{d_i} = \sum_{k=d_{min}}^{\infty} \bar{F}(x_n)^k u_k \\
&\leq  n\bar{F}(x_n)^{d_{min}} \frac{u_{d_{min}}}{n}+ \bar{F}(x_n) n\bar{F}(x_n)^{d_{min}}\frac{\sum_{k\geq0} u_{k}}{n} \\
\end{align*}
By Condition 1,
\begin{equation*}
 \frac{u_{d_{min}}}{n} \to p_{d_{min}} \text{ and } 
\frac{\sum_{k\geq0} u_{k}}{n} = \frac{\sum_i d_i^{(n)}}{n} \to \mathbb{E}[D] 
\end{equation*}
which, along with (\ref{xncon}) and (\ref{dicon}), shows that $S \to xp_{d_{min}}$. Thus, by Lemma~\ref{lemma1},
\begin{equation}
\label{dv2}
  d_{var}(P_{xp_{d_{min}}},P_{S}) \to 0
\end{equation}
Moreover by triangle inquality,
\begin{align*}
 d_{var}(\mathcal{L}(N_n\big|RCM[n]),P_{xp_{d_{min}}}) &\leq d_{var}(\mathcal{L}(N_n\big|RCM[n]),P_{S})\\
&\quad+d_{var}(P_{xp_{d_{min}}},P_{S})
\end{align*}
Therefore, from (\ref{dv1}) and (\ref{dv2}), we have finally that,
\begin{equation*}
  d_{var}(\mathcal{L}(N_n\big|RCM[n]),P_{xp_{d_{min}}})  \to 0
\end{equation*}
This completes the proof.

%
%
\proof[of corollary~\ref{mncor}]
It follows from the above theorem simply by taking 
\begin{equation}
\label{usex}
 x_n=\bar{F}^{-1}\Big[\big(\frac{x}{n}\big)^\frac{1}{d_{min}}\Big]
\end{equation}
and remarking that
\begin{equation*}
 \mathbb{P}(\text{\bf M}_n < x_n) = \mathbb{P}(N_n=0) \to e^{-xp_{d_{min}}}
\end{equation*}

\begin{rm1}
 For the case discussed in~\cite{Marc} where edge-weights are i.i.d exponential random variables with parameter $1$, i.e., $\bar{F}(x)=e^{-x}$, we have by Corollary~\ref{mncor} that for large $n$,
\begin{align*}
 &\mathbb{P}(ne^{-\text{\bf M}_n d_{min}} > O(\log n)) \approx e^{-O(\log n)p_{d_{min}}} \\
\text{or, }& \mathbb{P}(\text{\bf M}_n < \frac{1}{d_{min}}[\log n - o(\log n)]) \approx o(1) \\
\text{or, }& \mathbb{P}\big(\frac{\text{\bf M}_n}{\log n} \geq \frac{1}{d_{min}} - o(1)\big) \approx 1-o(1)
\end{align*}
which is consistent with the results in~\cite{Marc}. 
\end{rm1}

\begin{rm1}
 Elaborating on the conjecture regarding the longest edge of the MST on $\widetilde{RCM}[n]$, let $\widetilde{RCM}[n,x_n]$ denote the graph obtained by keeping only those edges of $\widetilde{RCM}[n]$ whose length is less than $x_n$.
The proof of Theorem~\ref{farout} suggests that for any $k \geq 2$, when $n \to \infty$ the connected components of size exactly $k$ diasppear in $\widetilde{RCM}[n,x_n]$. We confirm this below for $k=2$.
But there is nothing to suggest that there is one and only one unboundedly growing connected component, otherwise we could have concluded that far-out vertices are isolated in $\widetilde{RCM}[n,x_n]$ and thus
the longest edge of the MST on $\widetilde{RCM}[n]$ has the most far-out vertex as one of its ends, and its length is same as the nearest-neighbour distance of that vertex. 
\end{rm1} 

\begin{clus}
There are no clusters of size $2$ in $\widetilde{RCM}[n,x_n]$ w.h.p. 
\end{clus}
\proof
From the construction of $RCM[n]$, (\cite{Remco}):
\begin{equation*}
\mathbb{P}(i,j \text{ are connected in } RCM[n])=\frac{1}{cn}
\end{equation*}
where $c$ is some constant. Then, 
\begin{align*}
&\mathbb{P}(\exists \text{ a connected component of size } 2) \\
&\leq n^2\mathbb{P}(i,j \text{ are connected and they form an isolated component}) \\
&\leq n^2 \frac{1}{cn}F(x_n)\bar{F}(x_n)^{2d_{min}-1} \\
&\leq n\bar{F}(x_n)^{d_{min}}\bar{F}(x_n)^{2}
\end{align*}
which converges to $0$ by (\ref{xncon}) and (\ref{dicon}).

\section{Impact of clustering}
Taking cue from~\cite{perc-dcx} where (directionally) convex order as a measure of clustering provides a useful tool to compare point processes, we introduce a similar notion for the weighted $RCM$'s.
Given two random variables $X$ and $Y$, we say that $X$ is convexly smaller than $Y$, $X\leq_{cx} Y$,if $\mathbb{E}(f(X))\leq \mathbb{E}(f(Y))$ for all convex functions $f$. Remark that the comparison
necessitates $\mathbb{E}(X)= \mathbb{E}(Y)$. We will order $\widetilde{RCM}$'s according to the order of their limiting degree distributions. This order imposes further constraints on $\widetilde{RCM}$'s.
\begin{clus}
 Given two $\widetilde{RCM}$'s, $\widetilde{RCM}[\text{\bf d}^1]$ and $\widetilde{RCM}[\text{\bf d}^2]$, such that $D^1 \leq_{cx} D^2$, we should have one of the following two conditions:
\begin{enumerate}
 \item $d_{min}^1=d_{min}^2=d_{min}$ and $\mathbb{P}(D^1=d_{min})\leq \mathbb{P}(D^2=d_{min})$;
 \item $d_{min}^1>d_{min}^2$.
\end{enumerate}
\end{clus}
\proof
Suppose $d_{min}^1<d_{min}^2$. Then it is possible to construct a convex function $f$ such that $f(d_{min}^1)=1$ and $f(x)=0$ for all $x\geq d_{min}^1+1$. Then 
\begin{equation*}
\mathbb{P}(D^1=d_{min}^1)=\mathbb{E}(f(D^1))>\mathbb{E}(f(D^2))=0
\end{equation*}
which is a contradiction. Now for the case when $d_{min}^1=d_{min}^2=d_{min}$, we have
\begin{equation*}
 \mathbb{E}(z^{D^1-d_{min}}) \leq \mathbb{E}(z^{D^2-d_{min}})
\end{equation*}
for all $z\in (0,1)$. Letting $z$ approach $0$, we have $\mathbb{E}(z^{D^1-d_{min}}) \to \mathbb{P}(D^1-d_{min}=0)$ and therefore,
\begin{equation*}
 \mathbb{P}(D^1=d_{min})\leq \mathbb{P}(D^2=d_{min})
\end{equation*}

The above result along with Theorem~\ref{farout} suggests a connection between convex order on $\widetilde{RCM}$'s and the way they ``cluster''.

\section{Information transmission under channel-dependent emission capital constraints}
Consider $n$ servers labelled $1$ to $n$ such that at most $d_i$ channels of communication can pass through server $i$, where $(d_i)_1^n$ satisfy the conditions introduced earlier. A server is called \textit{closed} if the maximum possible channels of communication have been established through it, and is called \textit{available} otherwise.
  
Initially, a message is passed to a server $i$ from outside.
The server establishes channels with $d_i$ available servers  and then tries to transmit the message through each of them by spending an emission capital $C_n$ on each. But the transmission through channel $e_{ij}$ is successful only if $C_n$ exceeds $Y(e_{ij})$, where $\{Y(e_{ij})\}_{i,j}$ are i.i.d. random variables with distribution $F$. From then on, whenever the message
reaches server $k$ for the first time, it establishes channels with $d_k-1$ available servers (excluding the server from where it received the message) and tries to transmit the message through these channels in the same manner as described earlier. 

Our results show that for large $n$, $C_n$ should be at least $x_n$ as given by (\ref{usex}) if the message is to reach the entire network w.h.p.

\bibliographystyle{abbrv}

\end{document}